\documentclass[12pt,a4paper]{amsart}

\usepackage{color}
\usepackage{amscd,amssymb,amsopn,amsmath,amsthm,
                   graphics,amsfonts,enumerate,verbatim,calc,mathdots,mathrsfs}
\usepackage[dvips]{graphicx}
\usepackage{tikz}
\usetikzlibrary{intersections, calc, arrows}

\newtheorem{Theorem}{Theorem}[section]
\newtheorem{Proposition}[Theorem]{Proposition}
\newtheorem{Lemma}[Theorem]{Lemma}
\newtheorem{Corollary}[Theorem]{Corollary}

\newtheorem{Example}[Theorem]{Example}
\newtheorem{Definition}[Theorem]{Definition}

\makeatletter
  
  \@addtoreset{equation}{section}
\makeatother

\newcommand{\gp}{\mathfrak{p}}
\newcommand{\gq}{\mathfrak{q}}
\newcommand{\gm}{\mathfrak{m}}
\newcommand{\gn}{\mathfrak{n}}
\newcommand{\ga}{\mathfrak{a}}
\newcommand{\gb}{\mathfrak{b}}

\newcommand{\gM}{\mathfrak{M}}

\newcommand{\gA}{\mathfrak{A}}

\newcommand{\zz}{\mathbb{Z}}
\newcommand{\nn}{\mathbb{N}}

\newcommand{\gr}[1]{{\mathscr G}(#1)}

\newcommand{\mult}[2]{{\mathrm e}_{#1}(#2)} %multiplicity
\newcommand{\rmult}[2]{{\mathrm e}^R_{#1}(#2)} %relative multiplicity
\newcommand{\lc}[3]{{\mathrm H}^{\,#1}_{#2}(#3)} %local cohomology
\newcommand{\length}[2]{\ell_{#1}(#2)}
\newcommand{\ass}[2]{{\rm Ass}_{#1}\,{#2}}

\newcommand{\Supp}[2]{{\rm Supp}_{#1}\,#2}
\newcommand{\assc}[2]{{\rm Ass}^\circ_{#1}\,{#2}}
\newcommand{\Suppc}[2]{{\rm Supp}^\circ_{#1}\,#2}
\newcommand{\Min}[2]{{\rm Min}_{#1}\,{#2}}

\newcommand{\dep}[2]{{\rm depth}_{#1\,}#2}

\newcommand{\ra}{\longrightarrow}

\begin{document}
\title[Depth sensitivity of Hilbert coefficients]{
Depth Sensitivity of Hilbert Coefficients}

\author[K. Nishida]{Koji Nishida}%\footnotetext{  }
\address{
Department of Mathematics and Informatics,
Graduate School of Science,
Chiba University,
Yayoi-cho 1-33, Inage-ku,
Chiba 263-8522, Japan
}
\email{nishida@math.s.chiba-u.ac.jp}

\maketitle

\begin{abstract}
The purpose of this paper is to explain about the depth sensitivity
of the Hilbert coefficients defined for finitely generated graded modules
over graded rings.
The main result generalize the well known fact that 
the Cohen-Macaulayness of graded modules can be characterized 
using their multiplicities.
\end{abstract}

\section{Introduction}
Throughout this paper,
$R = \oplus_{n \in \nn\,}[ R ]_n$
is a Noetherian $\nn$-graded ring and
$M = \oplus_{n \in \nn\,}[ M ]_n$
is a finitely generated $\nn$-graded $R$-module
with $\dim_R M = s$,
where $\nn$ denotes the set of natural numbers including $0$.
Moreover, we always assume that
$[ R ]_0$ is an Artinian local ring whose residue field is infinite and
$R$ is generated by $[ R ]_1$ over $[ R ]_0$.
Then, as is well known, the integers
$\mult{0}{M}, \mult{1}{M}, \dots, \mult{s}{M}$
so called the Hilbert coefficients of $M$ are uniquely determined,
for which the following equality holds for $n \gg 0$;
\begin{equation}\label{eq1.1}
\sum_{k = 0}^n \length{[R]_0}{[M]_k} 
\,=\,
\sum_{i = 0}^s (-1)^i\cdot\mult{i}{M}\cdot\binom{n + s - i}{s - i},
\end{equation}
where $\length{[R]_0}{[M]_k}$ denotes the length of
$[ M ]_k$ as an $[ R ]_0$-module (c.f. \cite[Definition 4.1.5]{BH}).

The purpose of this paper is to show that
the Hilbert coefficients of $M$ enjoy the depth sensitivity,
introducing a new class of ssop (subsystem of parameters) for $M$.
In this paper, an ssop $f_1, \dots, f_n$ for $M$ is said to be admissible
if $f_1, \dots, f_n$ are all elements of $[R]_1$ and
$\dep{R_P}{M_P} \geq n$ for any $P \in \Supp{R}{M / (f_1, \dots, f_n)M}$
which does not include $[R]_1$.
Our main result can be described as follows.

\begin{Theorem}\label{1.1}
Let $f_1, \dots, f_{s - i}$ be an admissible ssop for $M$, where $0 \leq i < s$.
Then
\begin{itemize}
\item[{\rm (1)}]
$\mult{i}{M} \leq \mult{i}{M / (f_1, \dots, f_{s - i})M}$
if $i$ is even,
\item[{\rm (2)}]
$\mult{i}{M} \geq \mult{i}{M / (f_1, \dots, f_{s - i})M}$
if $i$ is odd, and
\item[{\rm (3)}]
$\mult{i}{M} = \mult{i}{M / (f_1, \dots, f_{s - i})M}$ holds
if and only if $\dep{R}{M} \geq s - i$.
\end{itemize}
\end{Theorem}

If $f_1, \dots, f_s$ is an sop (system of parameters) for $M$
consisting of elements of $[R]_1$,
then it is admissible obviously
as the homogeneous maximal ideal of $R$ is
the unique element of $\Supp{R}{M / (f_1, \dots, f_s)M}$,
and hence, from Theorem \ref{1.1},
it follows that
$\mult{0}{M} \leq \mult{0}{M / (f_1, \dots, f_s)M}$,
and the equality holds if and only if $M$ is a Cohen-Macaulay $R$-module.

In Section 3, we will show that
an ssop $f_1, \dots, f_n$ for $M$ such that
$f_1, \dots, f_n$ are all in $[R]_1$ is admissible
if and only if $(f_1, \dots, f_n)R$ is generated by a
superficial sequence,
which is a sequence defined by the condition a little bit weaker than 
that of a regular sequence.
And this fact will be the key point in our proof of Theorem \ref{1.1}.

Although quite many papers
have been been written already
on Hilbert coefficients , in almost all of them, 
the authors have studied mainly the invariants
$\mult{0}{I}, \mult{1}{I}, \dots, \mult{d}{I}$
that are defined for an $\gm$-primary ideal $I$
of a $d$-dimensional local ring $(A, \gm)$
as integers such that the equality
\begin{equation}\label{eq1.2}
\length{A}{A / I^{n + 1}} 
\,=\,
\sum_{i = 0}^d (-1)^i\cdot\mult{i}{I}\cdot\binom{n + d - i}{d - i}
\end{equation}
holds for $n \gg 0$
(cf. \cite{GN, HM, Hu, Nar, Nor, Sa, RV}).
Let us notice that $\mult{i}{I}$ can be understood
as the Hilbert coefficient of the associated graded ring
$\gr{I} = \oplus_{n \in \nn\,} I^n / I^{n + 1}$
in the context of this paper.
Moreover, making the theory applicable for finitely generated graded modules
enables us to compute the Hilbert coefficients more flexibly.
Of course, Theorem \ref{1.1} can be applied to
$\mult{i}{I}$ for $i \in \nn$ with $0 \leq i < d$ by considering a sequence
$a_1, \dots, a_{d - i}$ of elements of $I$ such that
$(a_1, \dots, a_{d - i})A$ is a part of minimal reduction of $I$ and
the leading forms of $a_1, \dots, a_{d - i}$ in $\gr{I}$
form an admissible ssop for $\gr{I}$.

In the last section, we will give several examples of computation of 
Hilbert coefficients for concrete quotient rings of polynomial rings over fields, that illustrate the depth sensitivity.

\section{Hilbert Coefficients}

In this section, let us briefly recall the definition and basic property of
the Hilbert coefficients of $M$,
setting $d = \dim R$ and $s = \dim_R M$.
We consider $s$ to be $-\infty$ if $M = 0$.

Let $t$ be an indeterminate and set
\[
P_M = \sum_{n \in \nn} \length{[R]_0}{[M]_n}\cdot t^n \in \zz[[ t ]],
\]
which is the Poincar\'{e} series of $M$.
If an exact sequence
$0 \rightarrow M_r \rightarrow \cdots \rightarrow M_1 \rightarrow M_0 \rightarrow 0$
of finitely generated $\nn$-graded $R$-modules is given,
it can be easily seen that the equality
$P_{M_0} - P_{M_1} + \cdots + (-1)^r\cdot P_{M_r} = 0$ holds.
On the other hand, we have
$P_{M(-r)} = t^r\cdot P_M$ for any $r \in \nn$,
where $M(-r)$ denotes the graded $\nn$-module such that
$[M(-r)]_n = M_{n - r}$ if $n \geq r$ and $[M(-r)]_n = 0$ if $0 \leq n < r$.
Applying these facts to the exact sequence
\[
0 \ra (0 :_M f)(-1) \ra M(-1) \stackrel{f}{\ra} M \ra M / fM \ra 0
\hspace{2ex} (f \in [R]_1),
\]
we get $(1 - t)\cdot P_M = P_{M / fM} - t\cdot P_{(0\,:_M f)}$.
Now we set
\begin{equation}\label{eq2.1}
\delta = 1 - t
\hspace{2ex}\mbox{and}\hspace{2ex}
\varphi_M = \delta^{\max\{s,\,0\}}\cdot P_M.
\end{equation}
Let us notice that $\varphi_M = P_M \in \zz[ t ]$ if $s \leq 0$.
Moreover, if $s > 0$,
we can choose an element $f \in [ R ]_1$ so that 
$\dim_R M / fM = s - 1$ as the residue field of $[ R ]_0$ is infinite,
and then it follows that $\dim_R (0 :_M f) < s$.
Therefore, by induction on $s$, we see $\varphi_M \in \zz[ t ]$.

When a polynomial $\xi \in \zz[ t ]$ is given,
let us denote by $\xi^{(i)}$ the $i$-th derivative of $\xi$
with respect to $t$ for $i \in \nn$.
Moreover, let $\xi(1)$ be the integer we get by
substituting $t = 1$ to $\xi \in \zz[ t ]$.
Then the Hilbert coefficients of $M$
can be defined as follows.

\begin{Definition}\label{2.1}
For any $i \in \nn$,
we set $\mult{i}{M} = 1 / i! \cdot \varphi_M^{\,(i)}(1)$ and
call it the $i$-th Hilbert coefficient of $M$.
\end{Definition}

Then, considering the Taylor expansion of $\varphi_M$ around $t = 1$,
we get
\[
\varphi_M = \sum_{i \in \nn} (-1)^i\cdot\mult{i}{M}\cdot\delta^i
\]
and $\mult{i}{M} \in \zz$ for any $i \in \nn$, which means
\begin{equation}\label{eq2.2}
\deg \varphi_M = \sup \{ i \in \nn \mid \mult{i}{M} \neq 0 \}.
\end{equation}
Here we notice that $\delta$ is a unit in $\zz[[ t ]]$ such that
$\delta^{-1} = 1 + t + t^2 + \cdots$.
If $M \neq 0$, then
\[
\delta^{-1}\cdot P_M = \delta^{-1 - s}\cdot \varphi_M =
\sum_{i \in \nn} (-1)^i\cdot\mult{i}{M}\cdot\delta^{i - s - 1}.
\]
For a formal power series $F \in \zz[[ t ]]$,
let us denote its coefficient of the term of $t^n$ by $F_n$,
i.e., $F = F_0 + F_1t + \cdots + F_nt^n + \cdots$.
Then we have
\begin{equation}\label{eq2.3}
(\delta^{-1}\cdot P_M)_n 
\,=\,
\sum_{i \in \nn} (-1)^i\cdot\mult{i}{M}\cdot(\delta^{i - s - 1})_n
\hspace{1ex}\mbox{fo any}\hspace{1ex}
n \in \nn.
\end{equation}
By straightforward calculation, we get
\[
(\delta^{-1}\cdot P_M)_n
\,=\,
\sum_{k = 0}^n \length{[R]_0}{[M]_k}.
\]
On the other hand,
if $i \leq s$,
then $s - i \geq 0$ and 
$\delta^{i - s - 1} = (1 + t + t^2 + \cdots)^{s - i + 1}$, 
and so
\[
(\delta^{i - s - 1})_n
\,=\,
\binom{n + s - i}{s - i}
\hspace{1ex}\mbox{for any}\hspace{1ex}
n \in \nn.
\]
Moreover, if $s < i \leq \deg \varphi_M$,
then $\delta^{i - s - 1} \in \zz[ t ]$ and
$\deg \delta^{i - s - 1} = i - s - 1 < \deg \varphi_M - s$,
which means
\[
(\delta^{i - s - 1})_n = 0
\hspace{1ex}\mbox{for any}\hspace{1ex}
n \in \nn
\hspace{1ex}\mbox{with}\hspace{1ex}
n \geq \deg \varphi_M - s.
\]
Therefore, (\ref{eq2.3}) implies that
the equality of (\ref{eq1.1}) holds
for any $n \geq \deg \varphi_M - s$
since $\mult{i}{M} = 0$ if $i > \deg \varphi_M$ by (\ref{eq2.2}).
Although the equality of (\ref{eq1.1}) implies that
$\mult{0}{M} = \length{R}{M}$ if $s = 0$,
but it is obvious by the definition of the $0$-th Hilbert coefficient.

\vspace{0.6em}
Next, we clarify what happens on  the Hilbert coefficients
when we shift the grading of modules.
For that purpose, let us recall that the binomial coefficients are defined
for $m \in \nn$ and $n \in \nn \setminus \{ 0 \}$ by
\[
\binom{m}{n} = \frac{1}{n!}\cdot\prod_{i = 0}^{n - 1}(m - i),
\]
which is equal to $0$ if $m < n$.
Moreover, we define
\[
\binom{m}{0} = 1
\]
for any $m \in \nn$.
Then, it can be easily seen that the equality
\[
\binom{m}{n} = \binom{m - 1}{n } + \binom{m - 1}{n - 1}
\]
which is very well known in the case where $m > n > 0$
still holds even if $0 < m \leq n$.
\begin{Lemma}\label{2.2}
For any $i, r \in \nn$, we have
\[
\mult{i}{M(-r)}
=
\sum_{j = 0}^i \binom{r}{j}\hspace{-0.5ex}\cdot\mult{i - j}{M}.
\]
In particular,
$\mult{0}{M(-r)} = \mult{0}{M}$ and
$\mult{i}{M(-1)} = \mult{i }{M} + \mult{i - 1}{M}$
if $i \geq 1$.
\end{Lemma}

{\it Proof}.\,
Because $P_{M(-r)} = t^r\cdot P_M$,
we have $\varphi_{M(-r)} = t^r\cdot \varphi_M$,
and so
\[
\mult{0}{M(-r)} = \varphi_{M(-r)}(1) =
1\cdot \varphi_M(1) = \mult{0}{M}.
\]
Hence the assertion holds certainly if $i = 0$.
Now we prove the assertion by induction on $r$
in the case where $i \geq 1$.
It is obvious if $r = 0$.
Moreover, as
\begin{eqnarray*}
\mult{i}{M(-1)} & = & \frac{1}{i!}\cdot(t\cdot \varphi_M)^{(i)}(1)
  \hspace{3ex}\mbox{by definition} \\
 & = & \frac{1}{i!}\cdot(t\cdot \varphi_M^{\,(i)} + i\cdot \varphi_M^{\,(i - 1)})(1)
 \hspace{3ex}\mbox{by Leibniz rule}  \\
 & = & \frac{1}{i!}\cdot 1 \cdot \varphi_M^{\,(i)}(1) +
                     \frac{1}{(i - 1)!}\cdot \varphi_M^{\,(i - 1)}(1) \\
 & = & \mult{i}{M} + \mult{i - 1}{M},
\end{eqnarray*}
the assertion holds if $r = 1$.
Let us assume $r \geq 2$ and set $N = M(-(r - 1))$.
Then
\[
\mult{i}{M(-r)} = \mult{i}{N(-1)} = \mult{i}{N} + \mult{i - 1}{N}.
\]
On the other hand, by the hypothesis of induction we have
\begin{eqnarray*}
\mult{i}{N} & = & \sum_{j = 0}^i
\binom{r - 1}{j}\cdot\mult{i - j}{M} \hspace{2ex}\mbox{and} \\
\mult{i - 1}{N} & = & \sum_{ k = 0}^{i - 1}
\binom{r - 1}{k}\cdot\mult{i - 1 - k}{M} =
\sum_{j = 1}^i
\binom{r - 1}{j - 1}\cdot\mult{i - j}{M}.
\end{eqnarray*}
Therefore it follows that
\[
\mult{i}{M(-r)} =  \mult{i}{M} + \sum_{j = 1}^i
\left\{\binom{r - 1}{j} + \binom{r - 1}{j - 1}\right\} \cdot\mult{i - j}{M}
= \sum_{j = 0}^i\binom{r}{j}\cdot\mult{i - j}{M}.
\]
Thus we have seen that the required assertion holds certainly.

\vspace{0.6em}
Here we introduce the Hilbert coefficients of another type.
Let us recall $d \geq \max\{s, 0\}$, and so $\delta^d\cdot P_M =
\delta^{d - \max\{s, 0\}}\cdot \varphi_M \in \zz[ t ]$.

\begin{Definition}\label{2.3}
For any $i \in \nn$,
we set $\rmult{\,i}{M} =
{1 / i!}\cdot(\delta^d\cdot P_M)^{(i)}(1)$ and call it
the $i$-th relative Hilbert coefficient of $M$.
\end{Definition}

If $s = d$,
$\mult{i}{M} = \rmult{i}{M}$ holds obviously for any $i \in \nn$.
In general, the relation between the Hilbert coefficients and 
the relative those can be described as follows.

\begin{Lemma}\label{2.4}
Suppose $M \neq 0$.
Then for any $i \in \nn$, we have
\[
\mult{i}{M} = (-1)^{d - s}\cdot\rmult{i + d - s}{M}
\hspace{2ex}\mbox{and}\hspace{2ex}
\rmult{\,i}{M} = \left\{\begin{array}{cl}
0 & \mbox{if $i < d - s$}, \vspace{0.5em}\\
(-1)^{d -s}\cdot\mult{i - d + s}{M} & \mbox{if $i \geq d - s$}.
\end{array}\right.
\]
\end{Lemma}

{\it Proof}.
Let us notice that
$\delta^d\cdot P_M = \delta^s\cdot\varphi_M$ as $s \geq 0$.
Then, by Leibniz rule, we have
\begin{eqnarray*}
(\delta^d\cdot P_M)^{(i)} & = & 
           \sum_{j = 0}^i \binom{i}{j}\cdot(\delta^{d - s})^{(j)}\cdot \varphi_M^{(i - j)} \\
 &  = & i!\,\cdot\hspace{-2.5ex}\sum_{j = 0}^{\min\{i,\,d - s\}} \hspace{-1.5ex}(-1)^j\cdot
           \delta^{d - s - j}\cdot\frac{1}{(i -j)!}\cdot \varphi_M^{(i - j)},
\end{eqnarray*}
and so, it follows that
\[
\rmult{\,i}{M} = \hspace{-2ex}\sum_{j = 0}^{\min\{i,\,d - s\}}\hspace{-2ex}(-1)^j\cdot
\delta^{d - s - j}(1)\cdot\mult{i - j}{M}.
\]
Hence we get the second equality as $\delta^{d - s - j}(1) = 0$ if $j < d - s$.
The first equality follows from the second one.

\vspace{0.6em}
The $i$-th relative Hilbert coefficients of
the $\nn$-graded $R$-modules appearing in an exact sequence
are related as follows.

\begin{Lemma}\label{2.5}
Let $0 \ra M_0 \ra M_1 \ra \cdots \ra M_n \ra 0$
be an exact sequence of finitely generated 
$\nn$-graded $R$-modules.
Then, for any $i \in \nn$, we have
\[
\sum_{k = 0}^n (-1)^k\cdot\rmult{\,i}{M_k} = 0.
\]
\end{Lemma}

{\it Proof}.
The  required equality holds since
\[
\sum_{k = 0}^n (-1)^k\cdot\frac{1}{\,i!\,}(\delta^d\cdot P_{M_k})^{(i)} = 
\left(\frac{\delta^d}{i!}\cdot\sum_{k = 0}^n(-1)^k\cdot
P_{M_k}\right)^{\hspace{-0.5ex}(i)} =
\left(\frac{\delta^d}{i!}\cdot 0\right)^{\hspace{-0.5ex}(i)} = 0.
\]

\vspace{0.5em}
The $i$-th relative Hilbert coefficient of the direct sum of
several $\nn$-graded $R$-modules is equal to the sum of
those of each modules as follows.

\begin{Lemma}\label{2.6}
Let $M_1, ... , M_n$ be finitely generated $\nn$-graded $R$-modules.
Then, for any $i \in \nn$, we have
\[
\rmult{\,i}{M_1 \oplus \cdots \oplus M_n} =
\rmult{\,i}{M_1} + \cdots + \rmult{\,i}{M_n}.
\]
\end{Lemma}

{\it Proof}.
If $n = 2$, we get the required equality 
applying Lemma \ref{2.5} to the exact sequence
$0 \ra M_1 \ra M_1 \oplus M_2 \ra M_2 \ra 0$.
Then the assertion can be verified easily by induction on $n$.

\vspace{0.6em}
The proof of the next result illustrates 
how to use the relative Hilbert coefficients.

\begin{Lemma}\label{2.7}
Let $f$ be an $M$-regular element of $[R]_k$,
where $0 < k \in \nn$.
Then, for any $i \in \nn$, we have
\[
\mult{i}{M / fM} = \sum_{j = 0}^i \binom{k}{j + 1}\cdot\mult{i - j}{M}.
\]
In particular,
$\mult{0}{M / fM} = k\cdot\mult{0}{M}$,
and $\mult{i}{M / fM} = \mult{i}{M}$
for any $i \in \nn$ if $k = 1$.
\end{Lemma}

{\it Proof}.
Because $d - \dim_R M / fM = d - s + 1$, by Lemma \ref{2.4}, 
we have 
\[
\mult{i}{M / fM} = (-1)^{d - s + 1}\cdot\rmult{j}{M / fM},
\]
where $j = i + d - s + 1$.
Moreover, applying Lemma \ref{2.5} to the exact sequence
\[
0 \ra M(-k) \stackrel{f}{\ra} M \ra M / fM \ra 0,
\]
we get $\rmult{j}{M / fM} = \rmult{j}{M} - \rmult{j}{M(-k)}$.
Consequently, it follows that
\begin{eqnarray*}
\mult{i}{M / fM} & = &
  (-1)^{d - s + 1}\cdot\{
  (-1)^{d - s}\cdot\mult{i + 1}{M} -
  (-1)^{d - s}\cdot\mult{i + 1}{M(-k)}\} \\
  & = & \mult{i + 1}{M(-k)} - \mult{i + 1}{M}
\end{eqnarray*}
as $j - d + s = i + 1$.
Therefore, by Lemma \ref{2.2}, we get
\[
\mult{i}{M / fM} = 
\sum_{j = 1}^{i + 1} \binom{k}{j}\cdot\mult{i + 1 - j}{M},
\]
whose right hand side is equal to that of the required equality.

\vspace{0.6em}
In the rest of this section,
we compute the Hilbert coefficients of 
concrete examples practically,
which is summarized as follows.

\begin{Example}\label{2.8}
Let $R$ be the polynomial ring over an infinite field $K$ 
with variables $x_1, \dots, x_d$.
Setting the degrees of $x_1, \dots, x_d$ are all $1$,
we regard $R$ as an $\nn$-graded ring such that $[R]_0 = K$.
\begin{itemize}
\item[{\rm (1)}]
$\displaystyle{
\mult{i}{R(-r)} = \binom{r}{i}
}$
for any $i, r \in \nn$.
\item[{\rm (2)}]
If $0 \neq f \in [R]_k$, where $0 < k \in \nn$,
then $\displaystyle{
\mult{i}{R / fR} = \binom{k}{i + 1}
}$
for any $i \in \nn$.
\vspace{0.3em}
\item[{\rm (3)}]
Let $f, g$ be an $R$-regular sequence such that 
$f \in [R]_k$ and $g \in [R]_\ell$ for some $0 < k, \ell \in \nn$.
Then, for any $i \in \nn$, we have
\begin{eqnarray*}
\mult{i}{M / (f, g)R} & = &
  \binom{k + \ell}{i + 2} - \binom{k}{i + 2} - \binom{\ell}{i + 2} \\
 & = & \sum_{j = 0}^i \binom{k}{i - j + 1}\cdot\binom{\ell}{j + 1}.
\end{eqnarray*}
\item[{\rm (4)}]
Let $\ga$ be the ideal of $R$ generated by the maximal minors of
an $m \times (m + 1)$ matrix whose entries are all linear forms of $R$.
Moreover, we assume ${\rm grade}\,\ga = 2$.
Then, for any $i \in \nn$, we have
\[
\mult{i}{R / \ga} = (i + 1)\cdot\binom{m + 1}{i + 2}.
\]
\end{itemize}
\end{Example}

{\it Proof}.
(1)\,\,
We put $\overline{R} = R / (x_1, \dots, x_d)R$.
Because $x_1, \dots, x_d$ is an $R$-regular sequence
and the degrees of $x_1, \dots, x_d$ are all $1$,
we have $\mult{i}{R} = \mult{i}{\overline{R}}$ for any $i \in \nn$
by Lemma \ref{2.7}.
On the other hand,
because $\overline{R}$ is concentrated in degree zero
and $[\overline{R}]_0 \cong K$,
we have $\varphi_{\overline{R}} = P_{\overline{R}} = 1$,
and hence 
\[
\mult{i}{\overline{R}} = \binom{0}{i} = \left\{\begin{array}{ccc}
1 & \mbox{if} & i = 0, \vspace{0.3em}\\ 
0 & \mbox{if} & i > 0.
\end{array}\right.
\]
Therefore the equality of (1) holds if $r = 0$.
Then we have
\begin{equation}\label{eq2.4}
\mult{i - j}{R} = \left\{\begin{array}{ccl}
1 & \mbox{if} & j = i \vspace{0.3em}\\
0 & \mbox{if} & 0 \leq j < i.
\end{array}\right.
\end{equation}
Hence, by Lemma \ref{2.2}, 
we see that the equality of (1) holds for any $r \in \nn$.

(2)\,\,
This assertion is deduced directly from Lemma \ref{2.7} and (\ref{eq2.4}).

(3)\,\,
Because $\dim R / (f, g)R = d - 2$,
we have
\begin{equation}\label{eq2.5}
\mult{i}{R / (f, g)R} = (-1)^2\cdot\rmult{i + 2}{R / (f, g)R}
\end{equation}
by Lemma \ref{2.4}.
Moreover, as $f, g$ is an $R$-regular sequence
such that $f \in [R]_k$ and $g \in [R]_\ell$,
$R / (f, g)R$ has an $\nn$-graded $R$-free resolution
of the following form;
\[
0 \ra R(-(k + \ell)) \ra R(-k) \oplus R(-\ell) \ra R
\ra R / (f, g)R \ra 0.
\]
Hence we see that the right hand side of (\ref{eq2.5}) is equal to
\[
\mult{i + 2}{R} - \mult{i + 2}{R(-k)} - \mult{i + 2}{R(-\ell)}
+ \mult{i + 2}{R(-(k + \ell))}
\]
by Lemma \ref{2.5}.
Consequently, by the formula given in (1) of this example, we get
\[
\mult{i}{R / (f, g)R} =
\binom{k + \ell}{i + 2} - \binom{k}{i + 2} - \binom{\ell}{i + 2}.
\]
On the other hand,
as $R / (f, g)R \cong (R / fR) / g(R / fR)$ and
$g$ is $R / fR$-regular, 
by Lemma \ref{2.7} and the formula given in (2) of this example, we get
\[
\mult{i}{R / (f, g)R} =
\sum_{j = 0}^i \binom{\ell}{j + 1}\cdot\mult{i - j}{R / fR} =
\sum_{j = 0}^i \binom{k}{i - j + 1}\cdot\binom{\ell}{j + 1}.
\]

(4)\,\,
Because $\dim R / \ga = d - 2$, we have
\begin{equation}\label{eq2.6}
\mult{i}{R / \ga} = (-1)^2\cdot\rmult{i + 2}{R / \ga}
\end{equation}
by Lemma \ref{2.4}.
Moreover, by the Burch's theorem, $R / \ga$ has an $\nn$-graded
$R$-free resolution of the following form;
\[
0 \ra R(-(m + 1))^{\oplus m} \ra R(-m)^{\oplus (m + 1)}
\ra R \ra R / \ga \ra 0.
\]
Hence we see that the right hand side of (\ref{eq2.6}) is equal to
\[
\mult{i + 2}{R} - (m + 1)\cdot\mult{i + 2}{R(-m)} + m\cdot\mult{i + 2}{R(-(m + 1))}
\]
by Lemma \ref{2.5}.
Consequently, by the formula given in (1) of this example, we get
\[
\mult{i}{R / \ga} = -(m + 1)\cdot\binom{m}{i + 2} + m\cdot\binom{m + 1}{i + 2},
\]
whose right hand side is equal to that of the equality of (4) 
as is verified by direct calculation.
Thus the proof is complete.

\section{Admissible ssop for graded module}

In this section, let us consider about a characterization
of admissible ssop for $M$,
which is necessary in the proof of Theorem \ref{1.1}.
We set $s = \dim_R M$ and
\[
\Suppc{R}{M} =
\{ P \in \Supp{R}{M} \mid [R]_1 \not\subseteq P \},
\hspace{2ex}
\assc{R}{M} = \{ P \in \ass{R}{M} \mid [R]_1 \not\subseteq P \}.
\]
Let us begin with the definition of superficial elements.

\begin{Definition} \label{3.1}{\rm (cf. \cite[Section 22]{Nag})}
An element $g \in [R]_1$ is said to be superficial for $M$ if
$(0 :_M g) \cap [M]_n = 0$ for $n \gg 0$.
\end{Definition}

It can be easily seen that the following three conditions;
$s \leq 0$,
$\Suppc{R}{M} = \phi$, and
$\assc{R}{M} = \phi$
are equivalent to each other,
and when this is the case,
any element of $[R]_1$ is superficial for $M$.
Even if $s > 0$,
the existence of superficial elements for $M$ is assured by 
the implication (2) $\Rightarrow$ (1) of the following lemma
as we are assuming that the residue field of $[R]_0$ is infinite.

\begin{Lemma}\label{3.2}
Suppose $s > 0$ and $g \in [R]_1$.
Then the following conditions{\rm ;}
\begin{itemize}
\item[{\rm (1)}]
$g$ is superficial for $M$,
\item[{\rm (2)}]
$g \not\in Q$ for any $Q \in \assc{R}{M}$,
\item[{\rm (3)}]
$g$ is $M_P$-regular for any $P \in \Suppc{R}{M / gM}$
\end{itemize}
are equivalent, and when this is the case,
we have $\dim_R M / gM = s - 1$.
\end{Lemma}

{\it Proof}.
Let us notice that $\assc{R}{M} \neq \phi$ as $s > 0$.
For $k \in \nn$, we set $M_{\geq k} = [M]_k \oplus [M]_{k + 1} \oplus \cdots$,
which is an $R$-submodule of $M$.

(1) $\Rightarrow$ (3) \,
Because $g$ is superficial for $M$,
we can choose $k$ so that $g$ is a non-zero-divisor on $M_{\geq k}$.
Now we take any $P \in \Suppc{R}{M / gM}$,
and choose $h \in [R]_1$ so that $h \not\in P$.
Then, as $h$ is a unit of $R_P$, we have
\[
M_P = h^k\cdot M_P = (h^k\cdot M)_P \subseteq (M_{\geq k})_P \subseteq M_P,
\]
and hence $M_P = (M_{\geq k})_P$.
Therefore $g$ is a non-zero-divisor on $M_P$,
which means that $g$ is $M_P$-regular as $g \in P$.

(3) $\Rightarrow$ (2)\,
Suppose that there exists $Q \in \assc{R}{M}$ such that $g \in Q$.
Then $QR_Q \in \ass{R_Q}{M_Q}$ and $Q \in \Suppc{R}{M / gM}$.
But these assertions contradicts to each other since
the second assertion implies that $g$ is $M_Q$-regular
by the condition (3).
Thus we see that the condition (2) is satisfied.
Then $\dim_R M / gM < s$,
and hence we get the last assertion of the lemma.

(2) $\Rightarrow$ (1)\,
Let $N_1 \cap \cdots \cap N_r = 0$ be a shortest primary decomposition
of $0$ in $M$.
We put $Q_i = \sqrt{N_i :_R M}$ for $i = 1, \dots, r$. 
Because the condition (2) implies that 
$g$ is $M$-regular if $\assc{R}{M} = \ass{R}{M}$, 
let us consider the case where $\assc{R}{M} \subsetneq \ass{R}{M}$.
Then $r \geq 2$, and by replacing the order of $N_1, \dots, N_r$ if necessary,
we may assume that
$\assc{R}{M} = \{ Q_1, \dots, Q_{r - 1} \}$.
Then $R_{\geq 1} \subseteq Q_r$.
Now we take an integer $m \gg 0$ so that
$Q_r^{\,m}\cdot M \subseteq N_r$.
Because $M$ is finitely generated as an $R$-module and
$R$ is generated by $[R]_1$ over $[R]_0$,
for some integer $k \gg 0$,
we have $M_{\geq k} \subseteq R_{\geq m}\cdot M =
(R_{\geq 1})^m\cdot M$,
which implies $M_{\geq k} \subseteq N_r$.
Here, let us take any $x \in (0 :_M g) \cap M_{\geq k}$.
Then, for any $i = 1, \dots, r - 1$,
we have $g \not\in Q_i$ and $gx = 0 \in N_i$,
which means $x \in N_i$ as $\ass{R}{M / N_i} = \{ Q_i \}$.
Hence we have
$x \in N_1 \cap \cdots \cap N_{r - 1} \cap N_r = 0$.
Therefore $(0 :_M g) \cap M_{\geq k} = 0$,
and $g$ is superficial for $M$ as required.

\begin{Definition}\label{3.3}
A sequence of elements $g_1, \dots, g_n$ in $[R]_1$
is said to be superficial for $M$ if
$g_i$ is superficial for $M / (g_1, \dots, g_{i - 1})M$
for any $i = 1, \dots, n$,
where $(g_1, \dots, g_{i - 1})M$ is considered to be $0$ if $i = 1$.
\end{Definition}

\begin{Lemma}\label{3.4}
Let $g_1, \dots, g_n$ be a superficial sequence for $M$,
where $0 < n \leq s$.
Then $\dim_R M / (g_1, \dots, g_n)M = s - n$ and
$g_1, \dots, g_n$ is an $M_P$-regular sequence for any
$P \in \Suppc{R}{M / (g_1, \dots, g_n)M}$.
\end{Lemma}

{\it Proof}.\,
We put $M_0 = M$ and  $M_i = M / (g_1, \dots, g_i)M$ for $i = 1, \dots, n$.
Then, as $g_i$ is superficial for $M_{i - 1}$ and
$M_{i - 1} / g_iM_{i - 1} \cong M_i$, we see that
$\dim_R M_i = \dim_R M_{i - 1} - 1$ and 
$g_i$ is $(M_{i - 1})_P$-regular for any 
$P \in \Suppc{R}{M_i}$
by Lemma \ref{3.2}.
Hence it follows that $\dim_R M_n = s - n$ and
$g_i$ is an regular element on $M_P / (g_1, \dots, g_{i - 1})M_P$
for any $P \in \Suppc{R}{M_n} \subseteq \Suppc{R}{M_i}$.
Thus we get the required assertion.

\vspace{0.6em}
We say that $M$ is a generalized Cohen-Macaulay $R$-module
if $\dep{R_P}{M_P} = s - \dim R / P$
for any $P \in \Suppc{R}{M}$
(cf. \cite{Tr2}).

\begin{Lemma}\label{3.8}
If $s > 0$ and $M$ is a generalized Cohen-Macaulay $R$-module,
then any ssop for $M$ consisting of elements of $[R]_1$ is a
superficial sequence for $M$.
\end{Lemma}

{\it Proof}.
Suppose that $0 < n \leq s$ and $f_1, \dots, f_n \in [R]_1$
is an ssop for $M$.
By induction on $n$, let us prove that
$f_1, \dots, f_n$ is a superficial sequence for $M$.
If $Q \in \assc{R}{M}$, we have
\[
0 = \dep{R_Q}{M_Q} = s - \dim R / Q,
\]
and hence $\dim R / Q = s$,
which means $f_1 \not\in Q$ as $f_1$ is an ssop for $M$.
Therefore, $f_1$ is a superficial element for $M$,
and hence the required assertion holds if $n = 1$.
So, let us consider the case where $n \geq 2$.
Then, by the hypothesis of induction,
it follows that $f_1, \dots, f_{n - 1}$ is a superficial sequence for $M$.
Here, let us take any $P \in \assc{R}{M / (f_1, \dots, f_{n - 1})M}$.
Then $\dep{R_P}{M_P / (f_1, \dots, f_{n - 1})M_P} = 0$ and
$f_1, \dots, f_{n - 1}$ is an $M_P$-regular sequence by Lemma \ref{3.4},
which means $\dep{R_P}{M_P} = n - 1$.
However, if $f_n \in P$ holds, we have
\[
\dep{R_P}{M_P} = s - \dim R / P \geq s - \dim_R M / (f_1, \dots, f_n)M
= s - (s - n) = n,
\]
which is impossible.
Thus we get $f_n \not\in P$,
which means that $f_n$ is a superficial element for $M / (f_1, \dots, f_{n - 1})M$.
Thus the proof is complete.

\begin{Definition}\label{3.5}
Let $s > 0$.
Then an ssop $f_1, \dots, f_n$ for $M$  such that
$f_i \in [R]_1$ for all $i$
is said to be admissible if
$\dep{R_P}{M_P} \geq n$ for any $P \in \Suppc{R}{M / (f_1, \dots, f_n)M}$.
\end{Definition}

\begin{Lemma}\label{3.6}
Let $s > 0$.
Then any sop for $M$ consisting of elements in $[R]_1$is admissible.
\end{Lemma}

{\it Proof}.\,
If $f_1, \dots, f_s$ is an sop for $M$,
the condition that
$\dep{R_P}{M_P} \geq s$ for any $P \in \Suppc{R}{M / (f_1, \dots, f_s)M}$
requires nothing since $\Suppc{R}{M / (f_1, \dots, f_s)M} = \phi$,
and hence the assertion holds.

\vspace{0.6em}
The next result guarantees the existence of an admissible ssop for $M$
consisting of $n$ elements for any integer $n$ with $0 < n \leq s$.
Moreover, it plays a key role in the proof of Theorem \ref{1.1}.

\begin{Proposition}\label{3.7}
Let $f_1, \dots, f_n$ be elements of $[R]_1$,
where $0 < n \leq s$.
Then the following conditions are equivalent.
\begin{itemize}
\item[{\rm (1)}]
$f_1, \dots, f_n$ is an admissible ssop for $M$.
\item[{\rm (2)}]
There exists a superficial sequence $g_1, \dots, g_n$ for $M$
generating $(f_1, \dots, f_n)R$.
\end{itemize}
\end{Proposition}

{\it Proof}.\,
The implication (2) $\Rightarrow$ (1) is a direct consequence of
Lemma \ref{3.4}. 
So let us prove (1) $\Rightarrow$ (2)
by induction on $\dim_R M$.

If $\dim_R M = 1$, then $n = 1$
and any $Q \in \assc{R}{M}$ does not contain
$f_1$ since $\dim R / Q = 1$ and $\dim_R M / f_1M = 0$,
which means that $f_1$ is a superficial element for $M$
by Lemma \ref{3.2}.
Hence the condition (2) is satisfied by taking $f_1$ itself as $g_1$.

So we assume that $\dim_R M = s \geq 2$ 
and the required implication is true
for finitely generated $\nn$-graded $R$-modules
of dimension $s - 1$.
Let us denote the maximal ideal of $[R]_0$ by $\gn$.
We set $\gA = (f_1, \dots, f_n)R$ and
$V = [\gA]_1 / \gn\cdot[\gA]_1$,
and we regard $V$ as a vector space over $[R]_0 / \gn$.
Because $f_1, \dots, f_n$ is an ssop for $M$ such that 
$f_i \in [R]_1$ for all $i$,
the dimension of $V$ as the vector space is $n$.
Moreover, for a homogeneous prime ideal $P$ of $R$,
we denote the image of $[P \cap \gA]_1$ under
the surjection $[\gA]_1 \ra V$ by $W_P$, 
which is a subspace of $V$.
Then, if $W_P = V$,
we have $[P \cap \gA]_1 = [\gA]_1$ by Nakayama's lemma,
which implies $\gA \subseteq P$.

Now, we take any $Q \in \assc{R}{M}$.
Then, as $\dep{R_Q}{M_Q} = 0 < n$,
we have $Q \not\in \Suppc{R}{M / \gA M}$ by the condition (1),
which means $W_Q \subsetneq V$ as $\gA \not\subseteq Q$.
Consequently, it follows that
we can choose an element $g_1 \in [\gA]_1 \setminus \gn\cdot[\gA]_1$
so that its image in $V$ is not contained in $W_Q$
for any $Q \in \assc{R}{M}$ as $[R]_0 / \gn$ is infinite.
Then we have $g_1 \not\in Q$ for any $Q \in \assc{R}{M}$,
and hence $g_1$ is a superficial element for $M$.

Let us choose elements $h_2, \dots, h_n$ of $[\gA]_1$
so that the images of $g_1, h_2, \dots, h_n$ in $V$ form a basis.
Then we have $\gA = (g_1, h_2, \dots, h_n)R$.
Here, we set $N = M / g_1M$.
It is obvious that $\dim_R N = s - 1$ and
$h_2, \dots, h_n$ is an ssop for $N$.
Moreover, if $Q \in \Suppc{R}{N / (h_2, \dots, h_n)N} =
\Suppc{R}{M / \gA M}$,
we have $\dep{R_Q}{N_Q} \geq n - 1$ since
$\dep{R_Q}{M_Q} \geq n$ by the condition (2) and
$g_1$ is $M_Q$-regular by Lemma \ref{3.2}.
Therefore, by the hypothesis of induction,
there exists a superficial sequence $g_2, \dots, g_n$ for $N$
such that $(h_2, \dots, h_n)R = (g_2, \dots, g_n)R$.
Then we see that $g_1, g_2, \dots, g_n$ is a superficial sequence for $M$
such that $\gA = (g_1, g_2, \dots, g_n)R$,
and the proof is complete.

\section{Proof of Theorem 1.1}

In this section, let us prove Theorem \ref{1.1}.
We set $d = \dim R, s = \dim_R M$ 
and $\gM$ to be the maximal homogeneous ideal of $R$.
First, we state a result which is a graded module version of Sally's machine
(cf. \cite[Lemma 2.2]{HM}, \cite[Lemma 1.4]{RV}).

\begin{Lemma}\label{4.1}
Let $g \in [R]_1$ be a superficial element for $M$.
We assume that $\dep{R}{M} > 0$ or $\dep{R}{M / gM} > 0$.
Then $g$ is $M$-regular,
and hence $\dep{R}{M} = \dep{R}{M / gM} + 1$.
\end{Lemma}

{\it Proof}.\,
If $\dep{R}{M} > 0$,
then $\assc{R}{M} = \ass{R}{M}$
as $\gM \not\in \ass{R}{M}$,
and hence $g$ is $M$-regular.
So we assume that $\dep{R}{M / gM} > 0$
in the rest of this proof.

Applying the local cohomology functor $\lc{0}{\gM}{\ast}$
to the exact sequence
$0 \ra gM \ra M \ra M / gM \ra 0$,
we get
\[
0 \ra \lc{0}{\gM}{gM} \ra \lc{0}{\gM}{M} \ra \lc{0}{\gM}{M / gM},
\]
which is an exact sequence of $\nn$-graded $R$-module of finite length.
Because $\dep{R}{M / gM} > 0$, we have $\lc{0}{\gM}{M / gM} = 0$,
and hence $\lc{0}{\gM}{gM} \cong \lc{0}{\gM}{M}$
as $\nn$-graded $R$-modules.
On the other hand,
applying the functor $\lc{0}{\gM}{\ast}$ to the exact sequence
$0 \ra (0 :_M g)(-1) \ra M(-1) \stackrel{g}{\ra} gM \ra 0$,
we get
\[
0 \ra \lc{0}{\gM}{0 :_M g}(-1) \ra \lc{0}{\gM}{M}(-1) \ra
\lc{0}{\gM}{gM} \ra \lc{1}{\gM}{0 :_M g}(-1),
\]
which is also an exact sequence of graded $R$-modules.
Here we notice that $0 :_M g$ has finite length as
$g$ is superficial.
Hence we have $\lc{0}{\gM}{0 :_M g} = 0 :_M g$ and
$\lc{1}{\gM}{0 :_M g} = 0$.
Thus we get an exact sequence
\[
0 \ra (0 :_M g)(-1) \ra \lc{0}{\gM}{M}(-1) \stackrel{g}{\ra} \lc{0}{\gM}{M}
\ra 0
\]
of $\nn$-graded $R$-modules of finite length.
It is obvious that the length does not change when we shift the grading
of any graded $R$-modules of finite length.
Consequently, it follows that
$\length{R}{0 :_M g} = 0$,
which means that $g$ is $M$-regular.
Thus the proof is complete.

\begin{Lemma}\label{4.2}
Suppose $0 < n < s$ and
$g_1, \dots, g_n \in [R]_1$
is a superficial sequence for $M$.
Then we have $\dep{R}{M} > n$ if and only if
$\dep{R}{M / (g_1, \dots, g_n)M} > 0$.
\end{Lemma}

{\it Proof}.\,
Let $M_0 = M$ and $M_j = M / (g_1, \dots, g_j)M$
for $j = 1, \dots, n$.
In order to get the required assertion,
it is enough to show that $\dep{R}{M_{j - 1}} > n - j + 1$
if and only if $\dep{R}{M_j} > n - j$ for $j = 1, \dots, n$,
which follows from Lemma \ref{4.1} since
$M_{j - 1} / g_jM_{j - 1} \cong M_j$ and
$g_j$ is a superficial element for $M_{j - 1}$.

\vspace{0.6em}
Next, we clarify how the Hilbert coefficients change 
when we take the quotient module by a superficial element.

\begin{Lemma}\label{4.3}
Suppose that $s > 0$ and
$g$ is a superficial element for $M$.
Then
\begin{itemize}
\item[{\rm (1)}]
$\mult{i}{M / gM} = \mult{i}{M}$ for any $i \in \nn$ with $i < s - 1$,
\item[{\rm (2)}]
$\mult{s - 1}{M / gM} = 
\mult{s - 1}{M} + (-1)^{s - 1}\cdot\length{R}{0 :_M g}$.
\end{itemize}
\end{Lemma}

{\it Proof}.\,
Let us take any $i \in \nn$.
Because $\dim_R M / gM = s - 1$ by Lemma \ref{3.2},
setting $j = d - (s - 1) = d - s + 1$, we have
\begin{equation}\label{eq4.1}
\mult{i}{M / gM} = (-1)^j\cdot\rmult{i + j}{M / gM}
\end{equation}
by Lemma \ref{2.4}.
Moreover, applying Lemma \ref{2.5} to the exact sequence
\[
0 \ra (0 :_M g)(-1) \ra M(-1) \stackrel{g}{\ra} M \ra M / gM \ra 0
\]
of $\nn$-graded $R$-modules,
we see that the right hand side of (\ref{eq4.1}) is equal to
\[
(-1)^j\cdot\{\rmult{i + j}{M} - \rmult{i + j}{M(-1)} +
\rmult{i + j}{(0 :_M g)(-1)}\}.
\]
Because $\dim_R M = \dim_R M(-1) = s$,
setting $k = d - s < i + j$, we have
\begin{eqnarray*}
 &   & (-1)^j\cdot\{\rmult{i + j}{M} - \rmult{i + j}{M(-1)}\}  \\
 & = & (-1)^j\cdot\{(-1)^k\cdot\mult{i + j - k}{M} -
   (-1)^k\cdot\mult{i + j -k}{M(-1)} \hspace{3ex}\mbox{by Lemma \ref{2.4}} \\
 & = & -\{\mult{i + 1}{M} - \mult{i + 1}{M(-1)}\}
           \hspace{3ex}\mbox{as $j + k$ is odd and $j - k = 1$} \\
 & = & \mult{i}{M} \hspace{3ex}\mbox{by Lemma \ref{2.2}.}
 \end{eqnarray*}
 Consequently, it follows that
 \[
 \mult{i}{M / gM} = \mult{i}{M} +
 (-1)^j\cdot\rmult{i + j}{(0 :_M g)(-1)}.
 \]
Let us notice that $\dim_R (0 :_M g) = 0$
as $g$ is a superficial element for $M$.
If $i < s - 1$,
we have $i + j < d = d - 0$,
and hence we get
$\rmult{i + j}{(0 :_M g)(-1)} = 0$
by Lemma \ref{2.4}.
Thus we see that the assertion (1) is true.
On the other hand, if $i = s - 1$,
we have $i + j = d$, and hence we get
\[
(-1)^j\cdot\rmult{i + j}{(0 :_M g)(-1)} =
(-1)^j\cdot(-1)^d\cdot\mult{0}{(0 :_M g)(-1)} =
(-1)^{j + d}\cdot\length{R}{0 :_M g}
\]
by Lemma \ref{2.4} and Lemma \ref{2.2}.
Therefore we see that the assertion (2) is also true as
$(-1)^{j + d} = (-1)^{s - 1}$,
and the proof is complete.

\vspace{0.6em}
Now, we are ready to prove Theorem \ref{1.1}.

\vspace{0.6em}
{\it Proof of Theorem \ref{1.1}}.\,
By Proposition \ref{3.7},
there exists a superficial sequence $g_1, \dots, g_{s - i}$
such that $(f_1, \dots, f_{s - i})R = (g_1, \dots, g_{s - i})R$.
Let $j$ be any integer such that $1 \leq j \leq s - i$.
We put $M_0 = M$ and $M_j = M / (g_1, \dots, g_j)M$.
Then we have $\dim_R M_{j - 1} = s - j + 1$,
$g_j$ is a superficial element for $M_{j - 1}$,
and $M_{j - 1} / g_jM_{j - 1} \cong M_j$.
Hence, if $j < s - i$,
we have $i < \dim_R M_{j - 1} - 1$,
and hence $\mult{i}{M_j} = \mult{i}{M_{j - 1}}$
by (1) of Lemma \ref{4.3}.
Thus we get
\[
\mult{i}{M} = \mult{i}{M_{s - i - 1}}.
\]
Moreover, by (2) of Lemma \ref{4.3},
we see that
\[
\mult{i}{M_{s - i}} =
\mult{i}{M_{s - i - 1}} + 
(-1)^i\cdot\length{R}{0 :_{M_{s - i - 1}} g_{s - i}}
\]
as $\dim_R M_{s - i - 1} = i + 1$.
Consequently, it follows that
\[
\mult{i}{M} = \mult{i}{M_{s - i}} +
(-1)^{i + 1}\cdot\length{R}{0 :_{M_{s - i - 1}} g_{s - i}},
\]
from which we get the assertions (1) and (2) of Theorem \ref{1.1}.
Moreover, we see that $\mult{i}{M} = \mult{i}{M_{s - i}}$ holds
if and only if $g_{s - i}$ is $(M_{s - i - 1})$-regular,
which is equivalent to $\dep{R}{M_{s - i - 1}} > 0$ by Lemma \ref{4.1}.
Therefore, we get the assertion (3) of Theorem \ref{1.1}
since $\dep{R}{M} \geq s - i$ holds if and only if
$\dep{R}{M_{s - i - 1}} > 0$ by Lemma \ref{4.2},
and the proof of Theorem \ref{1.1} is complete.

\begin{Corollary}\label{4.4}
Suppose that $M$ is a generalized Cohen-Macaulay $R$-module and
$f_1, \dots, f_{s - i}$ is any ssop for $M$ consisting of elements of $[R]_1$,
where $0 \leq i < s$.
Then the assertions {\rm (1)}, {\rm (2)} and {\rm (3)} of Theorem {\rm \ref{1.1}} hold.
\end{Corollary}

{\it Proof}.\,
Because $f_1, \dots, f_{s - i}$ is admissible for $M$ by Corollary \ref{3.8},
the required assertion follows from Theorem \ref{1.1} directly.

\section{Examples}

Throughout this section,
$K$ is an infinite field,
and $x_j, y_k$ are indeterminates over $K$
for $j, k \in \nn$.
We will compute the Hilbert coefficients of two kinds of
quotient rings $\overline{R}$ of 
polynomial rings $R$ setting the degrees of variables are all $1$.
Moreover, taking concrete admissible ssop $f_1, \dots, f_{s - i}$
for $\overline{R}$ from $[R]_1$,
we will compute $\mult{i}{\overline{R} / (f_1, \dots, f_{s - i})\overline{R}}$,
where $0 \leq i < s = \dim \overline{R}$.
The calculation results given in this section
back up the correctness of the three assertions of Theorem {\rm \ref{1.1}}.

\begin{Example}\label{5.1}
Suppose $s, d \in \nn$ and $0 < s < d$.
We set $R = K[x_1, \dots, x_d]$, $\gm = (x_1, \dots, x_d)R$,
$\gp = (x_1, \dots, x_{d - s})R$ and $\overline{R} = R / \gm\gp$.
Then the following assrtions hold.
\begin{itemize}
\item[{\rm (1)}]
$\dim \overline{R} = s$, $\dep{}{\overline{R}} = 0$ and
$\overline{R}$ is a generalized Cohen-Macaulay $R$-module.
\item[{\rm (2)}]
${\displaystyle
\mult{i}{\overline{R}} = 
\left\{\begin{array}{cl}
1 & \mbox{if} \hspace{1ex} i = 0, \vspace{0.2em}\\
0 & \mbox{if} \hspace{1ex} 0 < i < s, \vspace{0.2em}\\
(-1)^s\cdot(d - s) & \mbox{if} \hspace{1ex} i = s.
\end{array}\right.
}$
\item[{\rm (3)}]
If $0 \leq i < s$,
then $x_{d - s + i + 1}, \dots, x_d$
is a superficial sequence for $\overline{R}$
consisting of $s - i$ elements in $[R]_1$.
\item[{\rm (4)}]
${\displaystyle
\mult{i}{\overline{R} / (x_{d - s + i +1}, \dots, x_d)\overline{R}} = 
\left\{\begin{array}{cl}
d - s + 1 & \mbox{if} \hspace{1ex} i = 0, \vspace{0.2em}\\
(-1)^i\cdot(d - s) & \mbox{if} \hspace{1ex} 0 < i < s.
\end{array}\right.
}$
\end{itemize}
\end{Example}

{\it Proof}.\,(1)\,
Because $\Min{R}{\overline{R}} = \{ \gp \}$ and
\begin{equation}\label{eq5.1}
R / \gp \cong K[x_{d - s + 1}, \dots, x_d],
\end{equation}
we have $\dim \overline{R} = \dim_R \overline{R} =
\dim R / \gp = s$.
On the other hand,
as $\gp / \gm\gp \subseteq \overline{R}$,
we have $\ass{R}{\overline{R}} \supseteq \ass{R}{\gp / \gm\gp}
= \{ \gm \}$, and hence
$\dep{}{\overline{R}} = \dep{R}{\overline{R}} = 0$.
Moreover, if $P \in \Suppc{R}{\overline{R}}$,
we have $(\overline{R})_P \cong R_P / \gp R_P$
and $\gp R_P$ is generated by an $R_P$-regular sequence
consisting of $d - s$ elements,
and hence $\dep{R_P}{(\overline{R})_P} = \dim R_P - (d - s)
= (d - \dim R / P) - (d - s) = s - \dim R / P$.
Therefore $\overline{R}$ is a generalized Cohen-Macaulay $R$-module.

(2)\,
Because $\dim R - \dim \overline{R} = d - s$,
we have
\[
\mult{i}{\overline{R}} = (-1)^{d - s}\cdot\rmult{i + d -s}{\overline{R}}
\]
for any $i \in \nn$ by Lemma \ref{2.4}.
Moreover, applying Lemma \ref{2.5} to the exact sequence
\[
0 \ra \gp / \gm\gp \ra \overline{R} \ra R / \gp \ra 0
\]
of $\nn$-graded $R$-modules, we get
\[
\rmult{i + d - s}{\overline{R}} =
\rmult{i + d - s}{R / \gp} + \rmult{i + d - s}{\gp / \gm\gp}.
\]
Because we have the isomorphism of (\ref{eq5.1}),
by Lemma \ref{2.4} and Example \ref{2.8},
we get
\[
\rmult{i + d - s}{R / \gp} =
(-1)^{d - s}\cdot\mult{i}{R / \gp} =
\left\{\begin{array}{cl}
(-1)^{d - s} & \mbox{if} \hspace{1ex} i = 0, \vspace{0.2em}\\
 0 & \mbox{if} \hspace{1ex} i > 0.
\end{array}\right.
\]
On the other hand, if $i < s$,
we have $i + d - s < d = \dim R - \dim_R \gp / \gm\gp$,
and hence
\[
\rmult{i + d - s}{\gp / \gm\gp} = 0
\]
by Lemma \ref{2.4}.
Moreover, we have
\[
\rmult{d}{\gp / \gm\gp} = (-1)^d\cdot\mult{0}{\gp / \gm\gp}
= (-1)^d\cdot\length{R}{\gp / \gm\gp} = (-1)^d\cdot (d - s).
\]
Consequently, it follows that
\[
\mult{i}{\overline{R}} = \left\{\begin{array}{cl}
(-1)^{d - s}\cdot (-1)^{d - s} & \mbox{if} \hspace{1ex} i = 0, \vspace{0.2em}\\
0 & \mbox{if} \hspace{1ex} 0 < i < s, \vspace{0.2em}\\
(-1)^{d - s}\cdot (-1)^d\cdot (d - s) & \mbox{if} \hspace{1ex} i = s,
\end{array}\right.
\]
and hence we get the equalities of (2).

(3)\,
Because $\gm\gp + (x_{d - s + 1}, \dots, x_d)R =
\gp^2 + (x_{d - s + 1}, \dots, x_d)R$ is $\gm$-primary,
we get $\dim \overline{R} / (x_{d - s + 1}, \dots, x_d)\overline{R} = 0$,
and hence $x_{d - s + 1}, \dots, x_d$ is an sop for $\overline{R}$.
Then the required assertion follows from Lemma \ref{3.8}.

(4)\,
Let us take any $i \in \nn$ with $0 \leq i < s$.
We set $R' = K[x_1, \dots, x_{d - s + i}]$,
$\gm' = (x_1, \dots, x_{d - s + i})R'$ and $\gp' = (x_1, \dots, x_{d - s})R'$.
Then we have
\[
\overline{R} / (x_{d - s + i + 1}, \dots, x_d)\overline{R} \,\cong\,
R / (\gm\gp + (x_{d - s + i + 1}, \dots, x_d)R) \,\cong\,
R' / \gm'\gp'.
\]
If $i = 0$, we have $R' = K[x_1, \dots, x_{d - s}]$ and
$\gm' = \gp' = (x_1, \dots, x_{d - s})R'$, and hence
\[
\mult{0}{R' / \gm'\gp'} = \length{R'}{R' / (\gm')^2} = d - s + 1.
\]
Let us consider the case where $i > 0$.
In this case, we have $d - s = (d - s + i) - i$ and
$0 < i < d - s + i$, and hence, by (2), we get
\[
\mult{i}{R' / \gm'\gp'} = (-1)^i\cdot\{(d - s + i) - i\} =
(-1)^i\cdot( d - s).
\]
Thus we have seen that the required equalities hold,
and the proof is complete.

\begin{Example}\label{5.2}
Let $r, s$ be integers such that $0 < r < s$.
We set $R = K[x_1, \dots, x_s, y_1, \dots, y_r]$,
$\gp = (x_1, \dots, x_s)R$,
$\gq = (y_1, \dots, y_r)R$,
$\overline{R} = R / \gp\gq$ and
$z_j = y_j - x_j$ for $j = 1, \dots, r$.
Then the following assertions hold.
\begin{itemize}
\item[{\rm (1)}]
$\dim \overline{R} = s$ and $\dep{}{\overline{R}} = 1$.
\item[{\rm (2)}]
$\displaystyle{
\mult{i}{\overline{R}} =
\left\{\begin{array}{cl}
1 & \mbox{if}\hspace{1ex} i = 0, \vspace{0.2em}\\
(-1)^{s - r} & \mbox{if}\hspace{1ex} i = s - r, \vspace{0.2em}\\
(-1)^{s + 1} & \mbox{if}\hspace{1ex} i = s, \vspace{0.2em}\\
0 & \mbox{if}\hspace{1ex} i \not\in \{0, s - r, s\}.
\end{array}\right.
}$
\vspace{0.3em}
\item[{\rm (3)}]
If $0 \leq i < s - r$,
then $x_{r + i + 1}, \dots, x_s, z_1, \dots, z_r$
is an admissible ssop for $\overline{R}$ consisting of
$s - i$ elements in $[R]_1$ and
\vspace{-0.3em}
\[
\mult{i}{\overline{R} / (x_{r + i + 1}, \dots, x_s, z_1, \dots, z_r)\overline{R}}
= \left\{\begin{array}{cl}
r + 1 & \mbox{if}\hspace{1ex} i = 0, \vspace{0.2em}\\
(-1)^i\cdot r & \mbox{if}\hspace{1ex} 0 < i < s - r.
\end{array}\right.
\]
\item[{\rm (4)}]
If $s - r \leq i < s$,
then $z_{r - s + i + 1}, \dots, z_r$
is an admissible ssop for $\overline{R}$ consisting of
$s - i$ elements of $[R]_1$ and
\vspace{-0.3em}
\[
\mult{i}{\overline{R} / (z_{r - s + i + 1}, \dots, z_r)\overline{R}}
= \left\{\begin{array}{cl}
(-1)^{s - r}\cdot r & \mbox{if}\hspace{1ex} i = s - r, \vspace{0.2em}\\
(-1)^i\cdot (s - 1 - i) & \mbox{if}\hspace{1ex} s - r < i < s.
\end{array}\right.
\]
\end{itemize}
\end{Example}

{\it Proof.}\,
(1)\,
Because $\Min{R}{\overline{R}} = \{ \gp, \gq \}$,
\begin{equation}\label{eq5.2}
R / \gp \cong K[y_1, \dots, y_r]
\hspace{3ex}\mbox{and}\hspace{3ex}
R / \gq \cong K[x_1, \dots, x_s],
\end{equation}
we have $\dim \overline{R} = \dim_R \overline{R} =
\max \{ \dim R / \gp, \dim R / \gq \} = s$.
On the other hand,
setting $\gm = (x_1, \dots, x_s, y_1, \dots, y_r)R$,
we have the exact sequence
\begin{equation}\label{eq5.3}
0 \ra \overline{R} \ra R / \gp \oplus R / \gq \ra R / \gm \ra 0
\end{equation}
of $\nn$-graded $R$-modules since $\gp\gq = \gp \cap \gq$ and
$\gp + \gq = \gm$,
from which we get $\dep{}{\overline{R}} = 1$ by the depth lemma
since $\dim{}{R / \gm} = 0$ and $\dep{R}{(R / \gp \oplus R / \gq)} =
\min \{ r, s \} = r > 0$.

(2)\,
Because $\dim R - \dim \overline{R} = r$, we have
\[
\mult{i}{\overline{R}} = (-1)^r\cdot\rmult{i + r}{\overline{R}}
\]
for any $i \in \nn$ by Lemma \ref{2.4}.
Moreover, applying Lemma \ref{2.5} and Lemma \ref{2.6}
to the exact sequence (\ref{eq5.3}), we get
\[
\rmult{i + r}{\overline{R}} = \rmult{i + r}{R / \gp} +
\rmult{i + r}{R / \gq} - \rmult{i + r}{R / \gm}.
\]
Because we have the isomorphisms of (\ref{eq5.2}),
by Lemma \ref{2.4} and Example \ref{2.8}, we get
\[
\begin{array}{l}
\rmult{i + r}{R / \gp} = (-1)^r\cdot\mult{i}{R / \gp}
= \left\{\begin{array}{cl}
(-1)^r & \mbox{if} \hspace{1ex} i = 0, \vspace{0.2em} \\
0 & \mbox{if} \hspace{1ex} i > 0,
\end{array}\right.
\vspace{0.5em}
\\
\rmult{i + r}{R / \gq} = \left\{\begin{array}{cl}
0 & \mbox{if} \hspace{1ex} i < s - r, \vspace{0.2em}\\
(-1)^s\cdot\mult{i + r - s}{R / \gq} & \mbox{if} \hspace{1ex} i \geq  s - r,
\end{array}\right. 
\vspace{0.5em}
\\
\mult{i + r - s}{R / \gq} = \left\{\begin{array}{cl}
1 & \mbox{if} \hspace{1ex} i = s - r, \vspace{0.2em}\\
0 & \mbox{if} \hspace{1ex} i > s - r,
\end{array}\right. \vspace{0.5em}
\end{array}
\]
Moreover, because $\dim R - \dim R / \gm = r + s$, we have
\[
\rmult{i + r}{R / \gm} = \left\{\begin{array}{cl}
0 & \mbox{if} \hspace{1ex} i < s, \vspace{0.2em}\\
(-1)^{r + s}\cdot\mult{i - s}{R / \gm} & \mbox{if} \hspace{1ex} i \geq s,
\end{array}\right.
\hspace{2ex}
\mult{i - s}{R / \gm} = \left\{\begin{array}{cl}
1 & \mbox{if} \hspace{1ex} i = s, \vspace{0.2em} \\
0 & \mbox{if} \hspace{1ex} i > s.
\end{array}\right.
\]
Consequently, it follows that
\[
\mult{i}{\overline{R}} = \left\{\begin{array}{cl}
(-1)^r\cdot (-1)^r & \mbox{if} \hspace{1ex} i = 0, \vspace{0.2em} \\
(-1)^r\cdot (-1)^s & \mbox{if} \hspace{1ex} i = s - r, \vspace{0.2em} \\
(-1)^r\cdot(-1)\cdot(-1)^{r + s} & \mbox{if} \hspace{1ex} i = s, \vspace{0.2em} \\
0 & \mbox{if} \hspace{1ex} i \not\in \{ 0, s - r, s \}.
\end{array}\right.
\]
and hence we get the required equalities of (2).

(3)\,
Suppose $0 \leq i < s - r$.
Then $r + i + 1\leq s$.
We set $\ga = (x_{r + i + 1}, \dots, x_s, z_1, \dots, z_r)R$ and
take any $P \in \Suppc{R}{\overline{R} / \ga\overline{R}}$.
Then $P \supseteq \gp\gq$.
If $P \supseteq \gp$,
then $y_j = x_j + z_j \in P$ for any $j = 1, \dots, r$,
and hence $P \supseteq \gp + \gq = \gm$,
which is impossible.
Consequently, it follows that $P \not\supseteq \gp$ and $P \supseteq \gq$,
and hence we have
$(\overline{R})_P \cong R_P / \gq R_P = R_P / (y_1, \dots, y_r)R_P$.
Then, as $y_1, \dots, y_r, x_{r + i + 1}, \dots, x_s, z_1, \dots, z_r$
is a part of a $K$-basis of $[R]_1$,
we see that $x_{r + i + 1}, \dots, x_s, z_1, \dots, z_r$
is an $(\overline{R})_P$-regular sequence,
and hence $\dep{R_P}{(\overline{R})_P} \geq s - i$.
Thus we see that $x_{r + i + 1}, \dots, x_s, z_1, \dots, z_r$
is an admissible ssop for $\overline{R}$
consisting of $s - i$ elements in $[R]_1$.

Let us compute $\mult{i}{\overline{R} / \ga\overline{R}}$.
Because $x_jy_k \equiv x_jx_k$ mod $z_kR$
for any $j = 1, \dots, s$ and $k = 1, \dots, r$, we have
$\gp\gq + \ga = (x_1, \dots, x_{r + i})R\cdot (x_1, \dots, x_r)R + \ga$,
from which we get
\[
\overline{R} / \ga\overline{R} \cong
R / ((x_1, \dots, x_{r + i})R\cdot (x_1, \dots, x_r)R + \ga).
\]
Here we set $R' = K[x_1, \dots, x_{r + i}]$,
$\gm' = (x_1, \dots, x_{r + i})R'$ and
$\gp' = (x_1, \dots, x_r)R'$.
Then, as $R = R'[x_{r + i + 1}, \dots, x_s, z_1, \dots, z_r]$, we have
\[
\overline{R} / \ga\overline{R} \cong R' / \gm'\gp',
\]
and hence it is enough to show
\[
\mult{i}{R' / \gm'\gp'} = \left\{\begin{array}{cl}
r + 1 & \mbox{if} \hspace{1ex} i = 0, \vspace{0.2em}\\
(-1)^i\cdot r & \mbox{if} \hspace{1ex} 0 < i < s - r.
\end{array}\right.
\]
If $i = 0$, then $R' = K[x_1, \dots, x_r]$ and
$\gp' = \gm' = (x_1, \dots, x_r)R'$,
and hence $\mult{0}{R' / \gm'\gp'} =
\length{R'}{R' / (\gm')^2} = r + 1$.
On the other hand, if $0 < i < s - r$,
the required equality follows from (2) of Example \ref{5.1}.

(4)\,
Suppose $s - r \leq i < s$.
Then $0 \leq r - s + i < r$.
We set $\gb = (z_{r - s + i + 1}, \dots, z_r)R$
and take any $P \in \Suppc{R}{\overline{R} / \gb\overline{R}}$.
Then $P$ includes $\gp$ or $\gq$.
If $P \supseteq \gq$,
then $P \not\supseteq \gp$ as $P \neq \gm = \gp + \gq$,
and hence $z_{r - s + i +1}, \dots, z_r$
is an $(\overline{R})_P$-regular sequence since
$(\overline{R})_P \cong R_P / (y_1, \dots, y_r)R_P$
and $y_1, \dots, y_r, z_{r - s + i + 1}, \dots, z_r$
is a part of a $K$-basis of $[R]_1$.
On the other hand, if $P \supseteq \gp$,
then $P \not\supseteq \gq$, and hence
$z_{r - s + i + 1}, \dots, z_r$ is again an
$(\overline{R})_P$-regular sequence since
$(\overline{R})_P \cong R_P / (x_1, \dots, x_s)R_P$ and
$x_1, \dots, x_s, z_{r - s + i + 1}, \dots, z_r$
is a part of a $K$-basis of $[R]_1$.
Thus we get $\dep{R_P}{(\overline{R})_P} \geq s - i$,
and hence $z_{r - s + i + 1}, \dots, z_r$
is an admissible ssop for $\overline{R}$ consisting of
$s - i$ elements in $[R]_1$.

Let us compute $\mult{i}{\overline{R} / \gb\overline{R}}$.

First, we consider the case where $i = s - r$.
In this case, we have $\gb = (z_1, \dots, z_r)R$ and
$\gp\gq + \gb = \gp\cdot (x_1, \dots, x_r)R + \gb$.
Here we set $R'' = K[x_1, \dots, x_s]$,
$\gm'' = (x_1, \dots, x_s)R''$ and $\gp'' = (x_1, \dots, x_r)R''$.
Then, as $x_1, \dots, x_s, z_1, \dots, z_r$ is a $K$-basis of $[R]_1$,
we have
\[
R = R''[z_1, \dots, z_r]
\hspace{3ex}\mbox{and}\hspace{3ex}
\overline{R} / \gb\overline{R} \cong
\frac{R}{(\gp\gq + \gb)} \cong
\frac{R / \gb}{(\gp\gq + \gb)/\gb} \cong
R'' / \gm''\gp'',
\]
and hence it follows that
$\mult{s - r}{\overline{R} / \gb\overline{R}} = (-1)^{s - r}\cdot r$
by (2) of Example \ref{5.1}.

Next, we assume $i > s - r$.
In this case, we have
$0 < r - s + i < r$ and
\[
\gp\gq + \gb = 
\gp\cdot(y_1, \dots, y_{r - s + i}, x_{r - s + i + 1}, \dots, x_r)R + \gb.
\]
Here we set $S = K[x_1, \dots, x_s, y_1, \dots, y_{r - s + i}]$ and
\[
\overline{S} =
S /(x_1, \dots, x_s) (y_1,  \dots, y_{r - s + i}, x_{r - s + i + 1}, \dots, x_r)S.
\]
Then we have $\overline{R} / \gb\overline{R} \cong \overline{S}$ since
$x_1, \dots, x_s, y_1,  \dots, y_{r - s + i}, z_{r - s + i + 1}, \dots, z_r$
is a $K$-basis of $[R]_1$.
Hence it is enough to show that
$\mult{i}{\overline{S}} = (-1)^i\cdot(s - 1 - i)$.
For that purpose, we compute the Poincar\'{e} series of $\overline{S}$.
For $n \in \nn$,
let $\mathcal{M}_n(y_1, \dots, y_{r - s + i})$ be the set of 
monomials in $y_1, \dots, y_{r - s + i}$ of degree $n$.
We define 
$\mathcal{M}_n(x_1, \dots, x_{r - s +i}, x_{r + 1}, \dots, x_s)$ similarly.
If $n \geq 2$, it can be easily seen that
\[
\mathcal{M}_n(y_1, \dots, y_{r - s + i}) \cup
\mathcal{M}_n(x_1, \dots, x_{r - s +i}, x_{r + 1}, \dots, x_s)
\]
coincides with the set of monomials in
$x_1, \dots, x_s, y_1, \dots, y_{r - s + i}$ of degree $n$
which are not included in
$[(x_1, \dots, x_s) (y_1,  \dots, y_{r - s + i}, x_{r - s + i + 1}, \dots, x_r)S]_n$,
and hence we get
\begin{eqnarray*}
\length{K}{[\,\overline{S}\,]_n} & = &
   \sharp\mathcal{M}_n(y_1, \dots, y_{r - s + i}) + 
   \sharp\mathcal{M}_n(x_1, \dots, x_{r - s +i}, x_{r + 1}, \dots, x_s) \\
 & = &
   \binom{n + r - s + i - 1}{r - s + i - 1} + \binom{n + i - 1}{i - 1}.
\end{eqnarray*}
Therefore we have
\[
P_{\overline{S}} = 1 + (r + i)t +
\sum_{n = 2}^\infty \left\{
\binom{n + r - s + i - 1}{r - s + i - 1} + \binom{n + i - 1}{i - 1}
\right\}t^n.
\]
Let us recall that
\[
\delta^{-k} = \sum_{n = 0}^\infty \binom{n + k - 1}{k - 1} t^n
\]
holds for any $k \in \nn$, where $\delta = 1 - t$.
Hence we have
\begin{eqnarray*}
P_{\overline{S}} & = & 1 + (r + i)t + \{\delta^{-(r - s  + i)} - 1 - (r - s + i)t\}
+ \{\delta^{-i} - 1 - it\} \\
 & = & \delta^{-i} + \delta^{s - r - i} - 1 + (s - i)t.
\end{eqnarray*}
Then, as $\dim \overline{S} = \dim \overline{R} / \gb\overline{R} = i$,
we get
\begin{eqnarray*}
\varphi_{\overline{S}} & = & \delta^i\cdot P_{\overline{S}} \\
  & = & 1 + \delta^{s - r} - \delta^i + (s - i)t\delta^i \\
  & = & 1 + (-1)^{s - r}\cdot(t - 1)^{s - r} + (-1)^i\cdot(s - 1 - i)\cdot(t - 1)^i \\
  &    &  \hspace{30ex}     +\,(-1)^i\cdot(s - i)\cdot(t - 1)^{i + 1}.
\end{eqnarray*}
Because $\mult{i}{\overline{S}}$ is equal to the coefficient of the term of
$(t - 1)^i$ in the Taylor expansion of $\varphi_{\overline{S}}$ around $t = 1$,
we see that $\mult{i}{\overline{S}} = (-1)^i\cdot(s - 1 - i)$,
and the proof is complete.

\vspace{0.6em}
The admissible ssop for $\overline{R}$ given in (3) of Example \ref{5.2}
is not a superficial sequence since $x_{r + i + 1}$ is included in
$\gp \in \assc{R}{\overline{R}}$.
Example \ref{5.2} in the case where $r = s$ is studied  in \cite{In} .

\end{document}